\input amstex
\input amsppt.sty

\def\phi{\varphi}
\def\epsilon{\varepsilon}
\def\Re{\operatorname{Re}}

\magnification=1200
\NoBlackBoxes
\document
\topmatter
\title
An example concerning Bergman completeness
\endtitle
\author
W\l odzimierz Zwonek
\endauthor
\address
Instytut Matematyki, Uniwersytet Jagiello\'nski, Reymonta 4,
30-059 Krak\'ow, Poland
\endaddress
\email zwonek\@im.uj.edu.pl
\endemail
\thanks
The was partially supported by the KBN grant No. 2 PO3A 017 14.
\endthanks

\abstract
We construct a bounded plane domain which is Bergman complete but for which
the Bergman kernel does not tend to infinity
as the point approaches the boundary.
\endabstract

\endtopmatter

The disc with center
at $a\in\Bbb C$ and
radius $r>0$ we denote by $\triangle(a,r)$. We denote also
$E:=\triangle(0,1)$.
For $a\in\Bbb C$,
$0<r<R\leq\infty$ we denote the annulus
$P(a,r,R):=\{z\in\Bbb C:r<|z-a|<R\}$.

Let $D$ be a
bounded domain in $\Bbb C^n$. Let us denote by
$L_h^2(D)$ square integrable holomorphic functions on $D$. $L_h^2(D)$
is a Hilbert space with the scalar product induced from $L^2(D)$.
Let us define the {\it Bergman kernel of $D$}
$$
K_D(z)=\sup\{\frac{|f(z)|^2}{||f||_{L^2(D)}^2}:f\in L^2_h(D),
f\not\equiv 0\}.
$$
For the basic properties of the Bergman kernel
and other functions introduced below see
e.g \cite{Jar-Pfl}.

It is well-known that $\log K_D$ is a smooth plurisubharmonic function.
Therefore, we may define
$$
\beta_D(z;X):=\left(\sum_{j,k=1}^n\frac{\partial^2\log K_D(z)}
{\partial z_j\partial \bar z_k}X_j\bar X_k\right)^{1/2},\;
z\in D,\;X\in\Bbb C^n.
$$
The function $\beta_D$ is a pseudometric called {\it the Bergman pseudometric}.

For $w,z\in D$ we put
$$
b_D(w,z):=
\inf\{L_{\beta_D}(\alpha)\},
$$
where the infimum is taken over piecewise $C^1$-curves
$\alpha:[0,1]\mapsto D$ joining $w$ and $z$ and
$L_{\beta_D}(\alpha):=\int_0^1\beta_D(\alpha(t);\alpha^{\prime}(t))dt$.

We call $b_D$ {\it the Bergman distance of $D$}.

A bounded domain $D$ is called
{\it Bergman complete } if any $b_D$-Cauchy sequence is convergent to
some point in $D$ with respect to the standard topology of $D$.

Any bounded Bergman complete domain is pseudoconvex.

The proof of the Bergman completeness is often based on the proof
of the convergence of the Bergman kernel to infinity as the point
approaches the boundary, i.e. the following property
$$
\lim_{D\owns z\to\partial D}K_D(z)=\infty.\tag{$\ast$}
$$
All known Bergman complete domains have the
property \thetag{$\ast$}. On the other hand there are domains satisfying
\thetag{$\ast$} which are not
Bergman complete (take the Hartogs triangle). Let us recall
some known results on Bergman completeness and the property
\thetag{$\ast$}:

-- if $D$ is a bounded hyperconvex domain in $\Bbb C^n$, then
$D$ satisfies \thetag{$\ast$} (see \cite{Ohs~2}) and $D$ is Bergman complete
(see \cite{B\l o-Pfl} and \cite{Her}),

-- if $D$ is a bounded domain in $\Bbb C$ satisfying \thetag{$\ast$}, then
$D$ is Bergman complete (see \cite{Chen~2}),

-- all other known examples of Bergman complete domains (i.e. non-hyperconvex)
satisfy \thetag{$\ast$}, too (see \cite{Chen~1}, \cite{Her},
\cite{Jar-Pfl-Zwo} and \cite{Zwo}).

As already mentioned it has not been
clear whether the condition \thetag{$\ast$}
is necessary for a domain to be Bergman complete. As we show below
it is not the case.
The example given by us is a bounded domain in $\Bbb C$
(Theorem 5). Let us
underline here that the domain is given completely effectively.
As a by-product
we also get an effective example of a bounded fat domain in $\Bbb C$
not satisfying
\thetag{$\ast$} (see Corollary 3). For the non-effective proof of
the existence of such a domain see \cite{Jar-Pfl-Zwo}.

Below we restrict our considerations only to one-dimensional domains.

For a domain $D\subset\Bbb C$ and a function $f\in\Cal O(D)$ we denote
$||f||_D:=||f||_{L^2_h(D)}$. For $f,g\in L_h^2(D)$ we denote
$\langle f,g\rangle_D:=\int_Df \bar gd\lambda_2$.

For a fixed point $z_0\in\Bbb C$, $0<r<\infty$ we define
$\Cal O_0(P(z_0,r,\infty))$ as the set of holomorphic functions
$\phi$ from $\Cal O(P(z_0,r,\infty))$ such that their Laurent expansion
in $P(z_0,r,\infty)$
is of the form $\phi(z)=\sum_{n=1}^{\infty}\frac{a_{-n}}{(z-z_0)^n}$.
For such a function we also denote $(\phi)_{-1}(z):=\frac{a_{-1}}{z-z_0}$
and $(\phi)_{-2}(z):=\sum_{n=2}^{\infty}\frac{a_{-n}}{(z-z_0)^n}$.

Let us formulate the following two simple estimates, which we shall use
very extensively in the sequel:

\proclaim{Lemma 1} Let $\phi\in\Cal O(\triangle(z_0,R))$ ($0<R<\infty$). Then
for any $0\leq r\leq R$ the following inequality holds:
$$
||\phi||_{\triangle(z_0,r)}^2\leq\frac{r^2}{R^2}||\phi||_{\triangle(z_0,R)}^2.
$$
Let $\phi\in\Cal O_0(P(z_0,r,\infty))$. Assume that
$r\leq s\leq t$ and $r<t$. Then the following inequality holds:
$$
||\phi||_{P(z_0,s,t)}^2\leq\frac{\log t-\log s}{\log t-\log r}
||\phi||_{P(z_0,r,t)}^2.
$$
\endproclaim

At this place let us write down some technical property that we shall use
in the sequel.
Namely, the function
$u(x):=\frac{\log x-\log b}{\log x-\log a}$, $x>1$,
where $0<a<b<1$, is increasing, so $u(x)\leq u(2)$, $x\in (1,2)$.
Moreover, $u(2)\leq 2 u(1)$, if $a$ and $b$ are small enough,
for instance if $a,b\leq\exp(-4)$.

Below we shall consider sequences of positive numbers
$0<r_j<s_j<t_j$, $j=1,2,\ldots$ and points $z_1,z_2,\ldots\in E$
such that $\bar\triangle(z_j,t_j)\cap \bar\triangle(z_k,t_k)=\emptyset$
for any $j,k=1,2,\ldots$,
$j\neq k$, $0\not\in \bar \triangle(z_j,r_j)$, $j=1,2,\ldots$. Additionally,
we assume
that $z_N\to 0$. Then for such a fixed system of sequences we define
domains
$$
D_N:=E\setminus(\bigcup_{j=N}^{\infty}\bar\triangle(z_j,r_j)\cup\{0\}),
\;N=1,2,\ldots.
$$
In the sequel we shall also denote $D:=D_1$.

Additionally, we make some assumptions of the purely technical character
that we impose on the sequences considered:
$$
t_j<\exp(-4),\;r_j^2<\frac{|z_j|^2}{2},\;\frac{r_j^2}{t_j^2}+\frac{s_j}{t_j}+
\sqrt{\frac{2\log s_j}{\log r_j}}<1,\;\frac{2\log t_j}{\log r_j}+
\sqrt{\frac{2\log s_j}{\log r_j}}+\frac{s_j}{t_j}<1,\;j=1,2,\ldots.\tag{1}
$$
Our first aim is to find some sufficient conditions for the system of
sequences considered above implying the following condition
$$
\liminf_{D\owns z\to 0}K_D(z)<\infty.\tag{2}
$$

\proclaim{Lemma 2} Assume the following inequalities:
$$
\gather
\sum_{N=1}^{\infty}\frac{s_N}{t_N}<\infty,\;
\sum_{N=1}^{\infty}\sqrt{\frac{\log s_N}{\log r_N}}<\infty,\tag{3}\\
\sum_{N=1}^{\infty}\frac{-1}{\log r_N}<\infty.\tag{4}
\endgather
$$
Then there is a positive constant $C$ such that
$$
K_D(z)\leq C(K_E(z)+\sum_{j=1}^{\infty}(\frac{1}{|z-z_j|^2(-\log r_j)}+
\frac{r_j^2}{(|z-z_j|^2-r_j^2)^2})),\; z\in D.
$$
\endproclaim
\proclaim{Corollary 3} Let $D$ be as above. Assume the convergence
as in \thetag{3} and \thetag{4}. Assume also that $z_N>0$, $N=1,2,\ldots$
and
$\sum_{N=1}^{\infty}(\frac{-1}{z_N^2\log r_N}+\frac{r_N^2}
{(z_N^2-r_N^2)^2})<\infty$.
Then \thetag{2} is satisfied.
\endproclaim
It is easy to see that having given a sequence $z_N\to 0$,
$0<z_N<1$, $N=1,2,\ldots$, one may
easily (completely effectively) construct sequence
$\{r_N\}$ such that
the assumptions from Corollary 3 are satisfied.
 
\demo{Proof of Corollary 3} In view of Lemma 2
for $-\frac{1}{2}<z<0$ the following inequalities hold:
$$
K_D(z)\leq C(K_E(-1/2)+\sum_{j=1}^{\infty}(\frac{-1}{\log r_jz_j^2}+
\frac{r_j^2}{(z_j^2-r_j^2)^2})).
$$
The last expression is finite by the assumption of the Corollary.
\qed
\enddemo
\demo{Proof of Lemma 2}
Fix for a while some $N>0$.
Consider arbitrary $F\in L_h^2(D_N)$. It is a simple consequence of
the Laurent expansion of $F$ in the annulus $P(z_N,r_N,t_N)$
that $F=f+g$ in $D_N$,
where $f\in\Cal O(D_{N+1})$ and
$g\in\Cal O_0(P(z_N,r_N,\infty))$. It is easy to see that
$f\in L_h^2(D_{N+1})$ and $g\in L_h^2(P(z_N,r_N,R))$,
where $1<R<\infty$.
In view of Lemma 1 we have
$$
||f||_{\triangle(z_N,r_N)}^2\leq\frac{r_N^2}{t_N^2}
||f||_{\triangle
(z_N,t_N)}^2\leq
\frac{r_N^2}{t_N^2}||f||_{D_{N+1}}^2.
$$

Consequently,
$$
||f||_{D_N}^2=||f||_{D_{N+1}}^2-||f||_{\triangle(z_N,r_N)}^2\geq
(1-\frac{r_N^2}{t_N^2})||f||_{D_{N+1}}^2.\tag{5}
$$
On the other hand Lemma 1 gives the following estimates
$$
\gather
||g||_{D_N}^2\geq||g||_{P(z_N,r_N,t_N)}^2=
||g||_{P(z_N,r_N,1+|z_N|)}^2
-||g||_{P(z_N,t_N,1+|z_N|)}^2\geq\\
(1-\frac{\log(1+|z_N|)-\log t_N}
{\log(1+|z_N|)-\log r_N})||g||_{P(z_N,r_N,1+|z_N|)}^2\geq
(1-2\frac{\log t_N}{\log r_N})||g||_{P(z_N,r_N,1+|z_N|)}^2.\tag{6}
\endgather
$$
Now we want to find some upper estimates for the scalar product.
$$
\multline
|\langle f,g\rangle_{D_N}|\leq|\langle f,g
\rangle_{P(z_N,r_N,s_N)}|+|\langle f,g\rangle_{D_N\setminus
\bar\triangle(z_N,s_N)}|
\leq\\
||f||_{P(z_N,r_N,s_N)}||g||_{P(z_N,r_N,s_N)}+
||f||_{D_N\setminus\bar
\triangle(z_N,s_N)}
||g||_{D_N\setminus\bar \triangle(z_N,s_N)}.
\endmultline
$$

Since
$$
||f||_{P(z_N,r_N,s_N)}^2\leq||f||_{\triangle(z_N,s_N)}^2\leq
\frac{s_N^2}{t_N^2}||f||_{\triangle(z_N,t_N)}^2\leq\frac{s_N^2}{t_N^2}
||f||_{D_{N+1}}^2
$$
and
$$
||g||_{D_N\setminus\bar \triangle(z_N,s_N)}^2
\leq||g||_{P(z_N,s_N,1+|z_N|)}^2\leq
\frac{\log(1+|z_N|)-\log s_N}{\log(1+|z_N|)-\log r_N}
||g||_{P(z_N,r_N,1+|z_N|)}^2,
$$
the following inequality holds
$$
\multline
|\langle f,g\rangle_{D_N}|\leq\\
\frac{s_N}{t_N}||f||_{D_{N+1}}
||g||_{P(z_N,r_N,1+|z_N|)}+\sqrt{\frac{2\log s_N}{\log r_N}}
||f||_{D_{N+1}}
||g||_{P(z_N,r_N,1+|z_N|)}\leq\\
\frac{1}{2}(\frac{s_N}{t_N}+\sqrt{\frac{2\log s_N}{\log r_N}})
(||f||_{D_{N+1}}^2+||g||_{P(z_N,r_N,1+|z_N|)}^2).
\endmultline
\tag{7}
$$

Since
$$
||F||_{D_N}^2=||f+g||_{D_N}^2=||f||_{D_N}^2+||g||_{D_N}^2+
2\Re\langle f,g\rangle_{D_N},
$$
the inequalities \thetag{5}, \thetag{6} and \thetag{7}
give the following estimates
$$
\multline
||F||_{D_N}^2\geq
||f||_{D_{N+1}}^2
(1-\frac{r_N^2}{t_N^2}-\frac{s_N}{t_N}
-\sqrt{\frac{2\log s_N}{\log r_N}})+\\
||g||_{P(z_N,r_N,1+|z_N|)}^2(1-\frac{2\log t_N}{\log r_N}-\frac{s_N}{t_N}-
\sqrt{\frac{2\log s_N}{\log r_N}}).
\endmultline
\tag{8}
$$

More generally, using the Laurent expansion of $F\in L_h^2(D)$ in any annulus
$P(z_j,r_j,t_j)$, $j=1,\ldots,N$  we may find
$F_j\in\Cal O_0(P(z_j,r_j,\infty))$ (the choice of this $F_j$ is independent
of $N$) and $F_0^N\in\Cal O(D_{N+1})$ such that
$F=F_0^N+F_1+\ldots+F_N$ on $D_1$,
$F-F_j$ extends to a function holomorphic on $D\cup\bar\triangle(z_j,r_j)$.
Note that $F_0^N\in L_h^2(D_{N+1})$ and
$F_j\in L_h^2(P(z_j,r_j,R))$, $r_j<R<\infty$, $j=1,\ldots,N$.

Then in view of the inequality obtained in \thetag{8} applied recursively
we get the following estimate
$$
\multline
||F||_D^2\geq\\
\sum_{k=1}^N(||F_k||_{P(z_k,r_k,1+|z_k|)}^2
(1-\frac{2\log t_k}{\log r_k}-\sqrt{\frac{2\log s_k}{\log r_k}}-
\frac{s_k}{t_k})
\prod_{j=1}^{k-1}(1-\frac{r_j^2}{t_j^2}-\frac{s_j}{t_j}-
\sqrt{\frac{2\log s_j}{\log r_j}}))+\\
||F_0^N||_{D_{N+1}}^2\prod_{j=1}^N(1-\frac{r_j^2}{t_j^2}-\frac{s_j}{t_j}-
\sqrt{\frac{2\log s_j}{\log r_j}})).
\endmultline
$$

The convergence of the series $\sum_{N=1}^{\infty}\frac{s_N}{t_N}$
implies the convergence of the series
$\sum_{N=1}^{\infty}\frac{r_N^2}{t_N^2}$ and, consequently,
\thetag{3} implies that the infinite product
$$
\prod_{j=1}^{\infty}(1-\frac{r_j^2}{t_j^2}-\frac{s_j}{t_j}-
\sqrt{\frac{2\log s_j}{\log r_j}})
$$
is positive.

Moreover, $\inf_{j=1,2,\ldots}\{1-\frac{2\log t_j}{\log r_j}-
\sqrt{\frac{2\log s_j}{\log r_j}}-\frac{s_j}{t_j}\}$ is positive.

This altogether gives the existence of an $\epsilon>0$ such that for any $N$
$$
||F||_{D}^2\geq\epsilon(||F_0^N||_{D_{N+1}}^2+
\sum_{j=1}^N||F_j||_{P(z_j,r_j,1+|z_j|)}^2).
\tag{9}
$$
Our next aim is to show the local convergence of $F_0^N$ to a function
$F_0$ holomorphic on $E_*=\bigcup_{N=1}^{\infty}D_N$. Then in view
of \thetag{9} this convergence will imply that $F_0\in L_h^2(E_*)$
(consequently, we may treat $F_0$ as an $L_h^2$-function on $E$).
Note that the desired convergence follows from the local uniform convergence
of the series $\sum_{j=k}^{\infty}|F_j(z)|$ on $D_k$ for any $k=1,2,\ldots$,
which is proven below.

When we prove the above convergence then
$F=F_0+\sum_{j=1}^{\infty}F_j$ on $D_1$ and the following estimate will hold:
$$
||F||_{D}^2\geq\epsilon(||F_0||_E^2+
\sum_{j=1}^{\infty}||F_j||_{P(z_j,r_j,1+|z_j|)}^2).
\tag{10}
$$

Let us introduce some auxiliary functions:
$$
\gather
\tilde k_{j,-1}(z):=\\
\sup\{\frac{|(\phi)_{-1}(z)|^2}
{||\phi_{-1}||_{P(z_j,r_j,1+|z_j|)}^2}:\;\phi\in\Cal O_0(P(z_j,r_j,\infty)),
\;(\phi)_{-1}\not\equiv 0\}=\\
\frac{1}{2\pi|z-z_j|^2(\log(1+|z_j|)-\log r_j)},\\
\tilde k_{j,-2}(z):=\\
\sup\{\frac{|(\phi)_{-2}(z)|^2}
{||(\phi)_{-2}||_{P(z_j,r_j,1+|z_j|)}^2}:\;
\phi\in\Cal O_0(P(z_j,r_j,\infty)),
\;(\phi)_{-2}\not\equiv 0\}
\endgather
$$
for any $j=1,2,\ldots$, $z\in P(z_j,r_j,1+|z_j|)$.

Simple computations of the $L^2$-norms (of the functions
$(\phi)_{-2}$ imply that
there is some constant $C$ (independent of $j$)
such that
$$
\multline
\tilde k_{j,-2}(z)\leq CK_{P(z_j,r_j,\infty)}(z)=\\
\frac{Cr_j^2}{|z-z_j|^4}K_E(\frac{r_j}{z-z_j})=\frac{Cr_j^2}{
\pi(|z-z_j|^2-r_j^2)^2},\;z\in P(z_j,r_j,1+|z_j|),\;j=1,2,\ldots.
\endmultline
$$
Note that for any $k\leq N$ the following inequalities hold
$$
\multline
(\sum_{j=k}^N|F_j(z)|)^2\leq(\sum_{j=k}^N(|(F_j)_{-1}(z)|+|(F_j)_{-2}(z)|))^2\leq\\
(\sum_{j=k}^N||(F_j)_{-1}||_{P(z_j,r_j,1+|z_j|)}\tilde k_{j,-1}^{1/2}(z)+
          ||(F_j)_{-2}||_{P(z_j,r_j,1+|z_j|)}
           \tilde k_{j,-2}^{1/2}(z)))^2\leq\\
           (\sum_{j=k}^N(\tilde k_{j,-1}(z)+
\tilde k_{j,-2}(z)))(\sum_{j=k}^N
(||(F_j)_{-1}||_{P(z_j,r_j,1+|z_j|)}^2+
          ||(F_j)_{-2}||_{P(z_j,r_j,1+|z_j|)}^2))=\\
          (\sum_{j=k}^N(\tilde k_{j,-1}(z)+
\tilde k_{j,-2}(z)))\sum_{j=k}^N
||F_j||_{P(z_j,r_j,1+|z_j|)}^2,\;z\in D_k,
          \endmultline
$$
which finishes the proof of the desired properties of $F_0$ and \thetag{10}
(use the estimates for $\tilde k_{j,-1},\tilde k_{j,-2}$, \thetag{9} and
use the condition \thetag{4} to get for any $k$ the local boundedness
of the last expression, independently of $N$).

Now we may prove the required estimate. In view of \thetag{10}
we get the following estimate
$$
\multline
K_D(z)=\sup
\{\frac{|F(z)|^2}{||F||_{D}^2}:F\in L_h^2(D),\;F\not\equiv 0
\}\leq
\frac{1}{\epsilon}\sup\\
\{\frac{(|F_0(z)|+
\sum_{j=1}^{\infty}(|(F_j)_{-1}(z)|+|(F_j)_{-2}(z)|))^2}
{||F_0||_E^2+\sum_{j=1}^{\infty}
(||(F_j)_{-1}||_{P(z_j,r_j,1+|z_j|)}^2+||(F_j)_{-2}||_{P(z_j,r_j,1+|z_j|)}^2)},
\;F\not\equiv 0\},\;z\in D
\endmultline
$$
(the functions $F_j$ in the formula above come from the decomposition of $F$
considered earlier). And then proceeding as earlier we have
$$
\multline
K_D(z)
\leq
\frac{1}{\epsilon}
\sup\\
     \{
      \frac
        {\left(||F_0||_EK_E^{1/2}(z)+
         \sum_{j=1}^{\infty}
          (||(F_j)_{-1}||_{P(z_j,r_j,1+|z_j|)}\tilde k_{j,-1}^{1/2}(z)+
          ||(F_j)_{-2}||_{P(z_j,r_j,1+|z_j|)}
           \tilde k_{j,-2}^{1/2}(z))\right)^2
        }
        {||F_0||_E^2+\sum_{j=1}^{\infty}
         (||(F_j)_{-1}||_{P(z_j,r_j,1+|z_j|)}^2+
           ||(F_j)_{-2}||_{P(z_j,r_j,1+|z_j|)}^2
        }
     \}\\
\leq
\frac{1}{\epsilon}(K_E(z)+
\sum_{j=1}^{\infty}(\tilde k_{j,-1}(z)+
\tilde k_{j,-2}(z))),\;z\in D,
\endmultline
$$
which finishes the proof of the lemma (use the estimates
for $\tilde k_{j,-1}$ and $\tilde k_{j,-2}$).
\qed
\enddemo

\subheading{Remark 4} Note that the technical assumptions
in \thetag{1} do not cause
loss of generality for $z$ from the neighbourhood of $0$
(in particular, it does cause any loss of generality in Corollary 3).
The convergence of the series in \thetag{3} and \thetag{4} implies that for
$j$ large enough the technical properties from \thetag{1}
are always satisfied.
Therefore, because of the localization principle of the Bergman kernel
(see \cite{Ohs~1})
the estimates as in Lemma 2 (for $z$ from the neighbourhood
of $0$) remain valid without these technical assumptions.

\bigskip

Let us formulate our main result.

\proclaim{Theorem 5} There is a bounded domain $D\subset\Bbb C$ such that
$\liminf_{z\to\partial D}K_D(z)<\infty$ and $D$ is Bergman complete.
\endproclaim
\demo{Proof} The domain stated in the theorem will be some of the
domains considered earlier defined as $D_1$. Below we shall
impose some conditions on the sequences implying
that the domain has the property as desired.
Certainly, the point at which the Bergman kernel will not tend to infinity
will be $0$ (all other points from the boundary force
the Bergman kernel to diverge to infinity while tending to them).

Let us start with a sequence $x_n:=\frac{1}{n^5}$, $n\geq 2$.
We also define
$n^5$ different points lying on the circle of radius $x_n$.
$$
z_{n,j}:=x_n\exp(i\frac{2j\pi}{n^5}),\;j=0,\ldots,n^5-1.
$$
Note that there is some $C>0$ such that
$\frac{1}{Cn^{10}}\leq|z_{n,k}-z_{n,j}|$ and $|z_{n,0}-z_{n,1}|\leq\frac{C}{n^{10}}$
for any $n$ and for any
$j,k=0,\ldots,n^5-1$, $j\neq k$. Define
$$
t_n:=\frac{1}{3Cn^{10}},\;r_n:=\exp(-n^{19}),\;s_n:=\exp(-n).
$$
We also define $y_n:=\frac{x_n+x_{n+1}}{2}$.

Note that for any $n$
$\bar \triangle(z_{n,j},t_n)\cap\bar \triangle(z_{n,k},t_n)=\emptyset$,
$j\neq k$.
We also easily see that for $n,m$ large enough (for $n,m\geq n_0\geq 2$)
the circles $\partial\triangle(0,y_m)$ are disjoint from the discs
$\bar\triangle(z_{n,j},r_n)$, $j=0,\ldots,n^5-1$.
Now we build a sequence $\{z_N\}$
by gluing together one by one the (finite) sequences
$\{z_{n,j}\}_{j=0}^{n^5-1}$ (starting with $n=n_0$). We associate to them the sequences
$t_n$, $r_n$ and $s_n$ in such a way that $r_N$ (respectively,
$s_N$, $t_N$), where $N$ is such that
$z_N$ is associated to $z_{n,j}$, equals $r_n$
(respectively, $s_n$, $t_n$).
For indices large enough the sequences satisfy the technical assumptions from
\thetag{1}.

Note that the convergence as in \thetag{3} and \thetag{4} for the sequences
just defined will be satisfied
when we prove that
$$
\sum_{n=n_0}^{\infty}n^5a_n<\infty,
$$
where $a_n$ equals $\frac{s_n}{t_n}$ or $\sqrt{\frac{\log s_n}{\log r_n}}$
or $\frac{-1}{\log r_n}$. One can easily verify that this is the case.

One may also check that
$$
\sup\{\sum_{n=n_0}^{\infty}(\frac{n^5}{|y_m-x_n|^2(-\log r_n)}+
\frac{n^5r_n^2}{(|y_m-x_n|^2-r_n^2)^2}):
m=n_0,n_0+1,\ldots\}<\infty.
$$
Therefore, applying Lemma 2,
we easily see that there is some $M_1<\infty$ such that
the following inequality holds
$$
K_D(z)<M_1\text{ for any $z\in\bigcup_{n=n_0}^{\infty}
\partial \triangle
(0,y_n)\subset D$}.\tag{11}
$$
It follows from \thetag{11} that the Bergman kernel of $D$ does not diverge
to infinity as the point approaches $0$.

On the other hand take any $n\geq n_0$ and take a point
$z\in D\cap\partial \triangle(0,x_n)$.

Then
$$
\multline
K_D(z)\geq\max\{
\frac{
     \frac{1}{|z-z_{n,j}|^2}
     }
     {
     ||\frac{1}{\cdot-z_{n,j}}||_D^2
     },\;j=0,\ldots,n^5-1\}\geq
\frac{
     \frac{n^{20}}{C^2}
     }
     {
     ||\frac{1}{\cdot-x_n}||_{P(x_n,r_n,2)}^2
     }
     \geq\\
\frac{n^{20}}{2\pi C^2(\log 2-\log r_n)}=
\frac{n^{20}}{2\pi C^2(n^{19}+\log 2)}\to_{n\to\infty}
\infty.
\endmultline\tag{12}
$$

Now we are ready to prove the Bergman completeness of $D$. Suppose that $D$ is
not Bergman complete. Then there is a Cauchy sequence $\{w_k\}$ with respect
to the Bergman distance converging to the boundary (in the natural topology).
It is easy to see that this sequence must converge to $0$.
Choosing if necessary a subsequence we get from
the definition of the Bergman distance that there are a constant $M_2<\infty$
and a continuous function $\gamma:[0,1)\to D$ such that
$\lim_{t\to 0}\gamma(t)=0$, $\gamma_{|[0,1-\epsilon]}$ is piecewise $C^1$
and $L_{\beta_D}(\gamma_{|[0,1-\epsilon]})<M_2$ for any $\epsilon\in (0,1)$.
Note that the graph of $\gamma$ must intersect any set
$\partial\triangle(0,x_n)\cap D$
for $n\geq n_1$ with some $n_1\geq n_0$. Denote this point
of intersection by $v_n$. Then it follows from the definition of the Bergman
distance that the sequence $\{v_n\}$ is a Cauchy sequence
with respect to the Bergman distance. But, additionally,
it follows from \thetag{12} that $K_D(v_n)$ tends to infinity
as $n$ goes to infinity.
We prove below that this is impossible, which will finish the proof.
We follow the ideas from \cite{Chen~2}.

By a result from \cite{Pfl} there is a function $f\in L_h^2(D)$
such that $||f||_D=1$ and $\frac{|f(v_{n_j})|^2}{K_D(v_{n_j})}\to 1$
for some subsequence $\{z_{v_j}\}$ (see \cite{Chen~1} or \cite{Chen~2}).
Since functions bounded near $0$ are
dense in $L_h^2(D)$ (see \cite{Chen~2}, Lemma 4),
there exists a function $g\in L_h^2(D)$ such that $||f-g||_D\leq\frac{1}{2}$
and $g$ is bounded near $0$.
Then we have
$$
\frac{1}{2}\geq||f-g||_D\geq\frac{|f(v_{n_j})-g(v_{n_j})|}
{\sqrt{K_D(v_{n_j})}}\geq
\frac{|f(v_{n_j})|}{\sqrt{K_D(v_{n_j})}}-
\frac{|g(v_{n_j})|}{\sqrt{K_D(v_{n_j})}}\to 1
$$
-- contradiction.
\qed
\enddemo

\subheading{Remark 6} The last part of the proof of Theorem 5 is based
on the density of functions from $L_h^2(D)$ locally bounded in $0$
in the space $L_h^2(D)$. We quoted in this context the result
from \cite{Chen~2}. Actually, we may prove this result
in the special case
of domains considered in Lemma 2 directly
with elementary methods (without the use of the solutions of the
$\bar\partial$-problem). Namely, it follows from considerations in the proof
of Lemma 2 that for any $F\in L_h^2(D)$ we have
$F=F_0+\sum_{j=1}^{\infty}F_j$,
where $F_j\in L_h^2(P(z_j,r_j,1+|z_j|))\cap
\Cal O_0(P(z_j,r_j,\infty))$, $j=1,2,\ldots$ are as in the proof
of Lemma 2, $F_0\in L_h^2(E)$
and the convergence of the series is locally uniform on $D$.
Define $G_N:=F_0+\sum_{j=1}^NF_j$. Then $G_N\to F$ locally uniformly.
Assume that
the convergence is in $L_h^2(D)$-norm. Then
$F_0$ may be approximated in $L^2_h(E)$ by bounded holomorphic functions
on $E$ and
the functions $F_j$ may be approximated in $L_h^2(P(z_j,r_j,1+|z_j|))$
by bounded holomorphic functions in $P(z_j,r_j,1+|z_j|)$, $j=1,2,\ldots$.
Consequently, this implies the density of $H^{\infty}(D)$ in $L_h^2(D)$
(so a little more than it follows from the result of Chen).
To finish the proof it is sufficient to show that
$\sum_{j=1}^{N}F_j$ tends to $\sum_{j=1}^{\infty}F_j$ in $L_h^2(D)$.
But this can be seen from the considerations similar to that in the proof
of Lemma 2. Namely take $1\leq k<l$. Then
$$
||F_k+F_{k+1}+\ldots+F_l||^2_{D_k}=
||F_k||_{D_k}^2+||F_{k+1}+\ldots+F_l||_{D_k}^2+
2\Re\langle F_k,F_{k+1}+\ldots+F_l\rangle_{D_k}.
$$
The last expression is not larger than (repeat the reasoning
from the proof of \thetag{7})
$$
(||F_k||_{P(z_k,r_k,1+|z_k|)}^2+||F_{k+1}+\ldots+F_l||_{D_{k+1}}^2)
(1+\frac{s_k}{t_k}+\sqrt{\frac{2\log s_k}{\log r_k}}).
$$
Repeating this reasoning we get that the last expression is not larger than
$$
\sum_{j=k}^l||F_j||_{P(z_j,r_j,1+|z_j|)}^2\prod_{m=k}^{\min\{l-1,j\}}
(1+\frac{s_m}{t_m}+\sqrt{\frac{2\log s_m}{\log r_m}}).
$$
The assumptions on the convergence from \thetag{3} and \thetag{9}
easily finish the proof.

\subheading{Acknowledgment} The author would like to thank Professors
Peter Pflug and Zbigniew B\l ocki for helpful discussions
on the subject of the paper.

\enddocument